\documentclass[journal]{IEEEtran}

\usepackage{graphicx} 
  \DeclareGraphicsExtensions{.eps}
  
  \usepackage{tikz}
  \usepackage{pgfplots,pgfplotstable}
  \usetikzlibrary{pgfplots.groupplots}
  \pgfplotsset{compat=1.5}

\usepackage{multirow}

\usepackage{algorithm}
\usepackage[noend]{algorithmic} 

\usepackage{subcaption} 

\usepackage[font=footnotesize]{caption} 

\usepackage{amsmath}

\usepackage{textcomp}
\usetikzlibrary{shapes,arrows}

\usepackage{empheq} 

\usepackage{mathrsfs} 
\newsavebox\foobox
\newlength{\foodim}

\usepackage{stfloats} 

\usepackage[hyphens]{url}
\usepackage[hidelinks]{hyperref}
\hypersetup{colorlinks=false,linkcolor=black,citecolor=black,filecolor=black,urlcolor=black,breaklinks=true}
\usepackage[normalem]{ulem} 

\usepackage{cite} 

\usepackage{color}

\newcommand{\norm}[1]{\left\lVert#1\right\rVert} 

\usepackage{cancel} 

\usepackage{array} 

\usepackage[export]{adjustbox} 

\begin{document}

\bstctlcite{IEEEexample:BSTcontrol}

\title{Conditions for Regional Frequency Stability in Power System Scheduling---Part II: Application to Unit Commitment}

%
\author{Luis~Badesa,~\IEEEmembership{Member,~IEEE,}
        Fei~Teng,~\IEEEmembership{Member,~IEEE,}
        and~Goran~Strbac,~\IEEEmembership{Member,~IEEE}
\thanks{A portion of this work has been supported by the UK Engineering and Physical Sciences Research Council under project `Integrated Development of Low-Carbon Energy Systems’ (IDLES), grant EP/R045518/1.}
\thanks{
The authors are with the Department of Electrical and Electronic Engineering, Imperial College London, SW7 2AZ London, U.K. (email: luis.badesa@imperial.ac.uk, f.teng@imperial.ac.uk, g.strbac@imperial.ac.uk).}
%
}

%
%

\markboth{IEEE Transactions on Power Systems, April~2021}%
{Shell \MakeLowercase{\textit{et al.}}: Bare Demo of IEEEtran.cls for IEEE Journals}
%



\maketitle

\begin{abstract}
In Part I of this paper we have introduced the closed-form conditions for guaranteeing regional frequency stability in a power system. Here we propose a methodology to represent these conditions in the form of linear constraints and demonstrate their applicability by implementing them in a generation-scheduling model. This model simultaneously optimises energy production and ancillary services for maintaining frequency stability in the event of a generation outage, by solving a frequency-secured Stochastic Unit Commitment (SUC). We consider the Great Britain system, characterised by two regions that create a non-uniform distribution of inertia: England in the South, where most of the load is located, and Scotland in the North, containing significant wind resources. Through several case studies, it is shown that inertia and frequency response cannot be considered as system-wide magnitudes in power systems that exhibit inter-area oscillations in frequency, as their location in a particular region is key to guarantee stability. In addition, securing against a medium-sized loss in the low-inertia region proves to cause significant wind curtailment, which could be alleviated through reinforced transmission corridors. In this context, the proposed constraints allow to find the optimal volume of ancillary services to be procured in each region.
\end{abstract}


%
\IEEEpeerreviewmaketitle

\section*{Nomenclature}
\addcontentsline{toc}{section}{Nomenclature}

\subsection*{Acronyms}
\begin{IEEEdescription}[\IEEEusemathlabelsep\IEEEsetlabelwidth{$\textrm{RoCoF}_{\textrm{m}}$}]
\setlength\itemsep{0.3em} 
\item[CCGT] Combined Cycle Gas Turbine.
\item[COI] Centre Of Inertia.
\item[GB] Great Britain.
\item[MILP] Mixed-Integer Linear Program.
\item[OCGT] Open Cycle Gas Turbine.
\item[PFR] Primary Frequency Response.
\item[RES] Renewable Energy Sources.
\item[RoCoF] Rate-of-Change-of-Frequency.
\item[SUC] Stochastic Unit Commitment.
\end{IEEEdescription}

\subsection*{Indices and Sets}
\begin{IEEEdescription}[\IEEEusemathlabelsep\IEEEsetlabelwidth{$\textrm{RoCoF}_{\textrm{m}}$}]
\setlength\itemsep{0.3em} 
\item[$g,\,\, \mathcal{G}$] Index, Set of generators.
\item[$i,\,\, j$] All-purpose indices.
\item[$n,\,\, \mathcal{N}$] Index, Set of nodes in the SUC scenario tree.
\vspace*{3mm}
\end{IEEEdescription}

\subsection*{Constants and Parameters}
\begin{IEEEdescription}[\IEEEusemathlabelsep\IEEEsetlabelwidth{$\textrm{RoCoF}_{\textrm{max}}$}]
\setlength\itemsep{0.3em} 
\item[$\Delta\tau(n)$] Time-step corresponding to node $n$ (h).
\item[$\Delta f_{\textrm{max}}$] Maximum admissible frequency deviation at the nadir (Hz).
\item[$\Delta f^{\textrm{ss}}_{\textrm{max}}$] Maximum admissible frequency deviation at quasi-steady-state (Hz).
\item[$\pi(n)$] Probability of reaching node $n$.
\item[$\omega_i$] Angular frequency of inter-area oscillations for region $i$ (rad/s).
\item[$\textrm{A}_i$] Amplitude of inter-area oscillations for region $i$ (Hz).
\item[$\textrm{c}^\textrm{LS}$] Cost of load shedding (\pounds/MWh).
\item[$\textrm{c}^{\textrm{m}}_g$] Marginal cost of generating units $g$ (\pounds/MWh).
\item[$\textrm{c}^{\textrm{nl}}_g$] No-load cost of generating units $g$ (\pounds/h).
\item[$\textrm{c}^{\textrm{st}}_g$] Start-up cost of generating units $g$ (\pounds).
\item[$\textrm{D}_i$] Load damping factor in region $i$ (\%/Hz).
\item[$\textrm{H}_g$] Inertia constant of generating units $g$ (s).
\item[$\textrm{H}_\textrm{L}$] Inertia constant of the outaged generator (s).
\item[$\textrm{m}_j$] Regression weights.
\item[$\textrm{P}^{\textrm{D}}_i$] Total demand in region $i$ (MW).
\item[$\textrm{P}_{\textrm{g}}^{\textrm{max}}$] Maximum power output of units $g$ (MW).
\item[$\textrm{P}_{\textrm{g}}^{\textrm{msg}}$] Minimum stable generation of units $g$ (MW).
\item[$\textrm{P}_{\textrm{L}}^{\textrm{max}}$] Rated power of the largest generator (MW).
\item[$\textrm{R}_g^\textrm{max}$] PFR capacity of generators $g$ (MW).
\item[$\textrm{R}_g^\textrm{slope}$] Proportion of headroom that can contribute to PFR.
\item[$\textrm{RoCoF}_{\textrm{max}}$] Maximum admissible RoCoF (Hz/s).
\item[$\textrm{V}_i$] Voltage magnitude in bus $i$ (kV).
\item[$\textrm{X}_{i,j}$] Reactance of the transmission line connecting buses $i$ and $j$ ($\Omega$).
\end{IEEEdescription}

\subsection*{Decision Variables}
\begin{IEEEdescription}[\IEEEusemathlabelsep\IEEEsetlabelwidth{$\textrm{RoCoF}_{\textrm{ma}}$}]
\setlength\itemsep{0.3em} 
\item[$N_g^\textrm{sg}(n)$] Number of units $g$ that start generating in node $n$.
\item[$N_g^{\textrm{up}}$] Number of online generating units of type $g$.
\item[$P_g$] Power produced by generating units $g$ (MW).
\item[$P^{\textrm{L}}$] Largest power infeed (MW).
\item[$P^\textrm{LS}$] Load shed (MW).
\item[$R_g$] PFR provision from generating units $g$ (MW).
\end{IEEEdescription}

\subsection*{Linear Expressions \normalfont{(linear combinations of decision variables)}}
\begin{IEEEdescription}[\IEEEusemathlabelsep\IEEEsetlabelwidth{$\textrm{RoCoF}_{\textrm{m}}$}]
\setlength\itemsep{0.3em} 
\item[$C_g(n)$] Operating cost of units $g$ at node $n$ in the SUC (\pounds).
\item[$H_i$] System inertia in region $i$ (MW$\cdot \textrm{s}$).
\item[$H$] Total system inertia (MW$\cdot \textrm{s}$).
\item[$R_i$] PFR from all providers in region $i$ (MW).
\item[$R$] Total system PFR (MW).
\end{IEEEdescription}

\subsection*{Functions and Operators}
\begin{IEEEdescription}[\IEEEusemathlabelsep\IEEEsetlabelwidth{$\textrm{RoCoF}_{\textrm{m}}$}]
\setlength\itemsep{0.3em} 
\item[$\norm{\;\cdot\;}$] $\ell^2\textrm{-norm}$.
\item[$\Delta f_i(t)$] Time-evolution of post-fault frequency deviation from nominal state in region $i$ (Hz).
\item[$t_\textrm{nadir}$] Time when the frequency nadir occurs (s).
\end{IEEEdescription}

\section{Introduction}
%
%
%
%
\IEEEPARstart{I}{n} order to assure a secure operation of the system from a frequency-performance point of view, system operators must procure certain frequency services that would only come into play in the event of a frequency drop.  

Frequency services are any type of ancillary service that helps comply with frequency regulation, that is, that helps restore a power equilibrium after a generation/demand outage. These services are: inertia (i.e.~the kinetic energy stored in the rotating masses of synchronous generators), 
frequency response (a power injection from different devices that is activated after a loss), 
load damping (typically provided by frequency-responsive loads that reduce their consumption after a frequency drop) and the size of the largest possible loss of generation or demand (corresponding to the \textit{N}-1 reliability criterion).

These ancillary services increment the operating cost of a power system, since the level of system inertia depends on the number of thermal units committed, while PFR is a function of the headroom in online plants. The main providers of frequency services are generators and loads,  
therefore these services are inherently related to energy provision. In order to optimise the provision of inertia and PFR, several works have focused on constraining Optimal Power Flow and Unit Commitment problems to explicitly respect frequency security. Different approaches have been proposed in the literature to limit the maximum admissible RoCoF, minimum admissible frequency at nadir and minimum admissible frequency at quasi-steady-state that would result from the optimisation solution. 

The authors in \cite{UCRestrepoGaliana} enforced a minimum threshold of PFR in a scheduling algorithm, a threshold driven by the requirement for frequency to return to its nominal value following a power outage (i.e.~the frequency steady-state limit). The RoCoF and frequency nadir limits were considered in \cite{OMalleyDeload}, by simulating a simplified model of the system dynamics.
An analytical constraint for the frequency nadir was introduced in an OPF in \cite{OPFChavez}.
While this formulation could only consider a fixed system inertia and optimise over the volume of PFR required, this limitation was overcome in \cite{FeiStochastic} using a linearisation technique on the originally non-convex constraint. 

Simultaneous optimisation of PFR and a fast frequency response service from batteries has been the recent focus. The authors in \cite{ERCOT_EFR} assume inertia to be given, and propose a linear replacement ratio to approximate the non-linear relationship between these two types of frequency response. Co-optimisation with inertia is achieved through analytical formulations for the frequency nadir in \cite{VincenzoEFR,LuisEFR}. The latest approaches \cite{LuisMultiFR,NewZealandMultiFR} allow to optimise multi-speed frequency response, not limited to the two speeds considered in earlier works.

Certain assumptions were made in these studies to overcome the mathematical difficulties and obtain relatively simple algebraic constraints suitable for implementation in convex optimisation problems. Analytical formulations rarely account for the support from load damping or consider the size of the largest loss as a decision variable, with the exception of \cite{LuisEFR}. All these studies assumed frequency to be equal throughout the grid.

In the present paper, a scheduling model is constrained to guarantee frequency stability in every region of the power system, therefore moving for the first time beyond the uniform frequency model considered in previous works. To highlight the implications of such model, several case studies are run considering the Great Britain power system, characterised by two regions that create a non-uniform distribution of inertia: England in the South, where most of the load is located, and Scotland in the North, containing significant wind resources. The results presented demonstrate that inertia and frequency response cannot be considered as system-wide magnitudes in power systems that exhibit inter-area oscillations in frequency, as their location in a particular region is key to guarantee stability. In this context, the proposed constraints allow to find the optimal volume of ancillary services to be procured in each region.

The contributions of Part II of this paper are the following:
\begin{enumerate}
    \item A numerical method that provides linear constraints for regional frequency stability is proposed, enabling to implement them as constraints in Mixed-Integer Linear Programs (MILPs). A regional-frequency-secured Unit Commitment is formulated for the first time.
    \item The consequences of inter-area oscillations in frequency are thoroughly studied using several relevant examples, presenting fundamental insight such as a quantification of the higher need for ancillary services as compared to the uniform frequency model. 
    \item The optimal location of inertia and frequency response is investigated, demonstrating that both the regional cost of each service and the size of the contingency in each region have a substantial impact on the optimal schedule. 
\end{enumerate}

The rest of this paper is organised as follows. Section~\ref{Sec:Applicability} demonstrates the applicability of the proposed frequency-stability conditions to be implemented in MILPs. Section~\ref{Sec:CaseStudies} presents the results of several relevant case studies, analysing the impact on power system scheduling of inter-areas oscillations in regional frequencies. 
Finally, Section~\ref{Sec:Conclusion} provides the conclusion and proposes future lines of work.

\section{Applicability of the constraints for regional frequency stability to scheduling problems} \label{Sec:Applicability}

In Part I of this paper we have deduced analytical constraints for guaranteeing that RoCoF and nadir stay within pre-defined limits in the event of any credible contingency. Here we propose a numerical method to formulate those constraints as linear expressions, a method consisting on using linear regressions on samples obtained from dynamic simulations of the power system. Although some conservative steps have been taken to deduce these linear constraints (i.e.~steps that tighten the feasible region for frequency security, therefore requiring a higher volume of frequency services than absolutely necessary), we also demonstrate that the resulting constraints represent the actual stability boundary accurately, showing only a small degree of conservativeness.

\subsection{RoCoF constraint: numerical estimation} \label{sec:RoCoFnumericalEstimation}

The analytical deduction of the RoCoF constraint for each region gives the following inequality:
\begin{equation} \label{eq:RoCoF_analytical_TwoRegions}
\left|\mbox{RoCoF}_i\right| = \frac{P^\textrm{L}}{2 (H_1+H_2)}+\textrm{A}_i\cdot\omega_i \leq \mbox{RoCoF}_{\textrm{max}}
\end{equation}

Here we propose a numerical-estimation method for the parameters of the inter-area oscillations appearing in eqs.~(\ref{eq:RoCoF_analytical_TwoRegions}), i.e.~parameters $\textrm{A}_i$ and $\omega_i$. This method allows to obtain linear RoCoF constraints for each region, in the form:

\begin{multline} \label{eq:RoCoF_regional_regression}
  \frac{P^\textrm{L}}{2 (H_1+H_2)} + \frac{\textrm{m}_1 H_1+\textrm{m}_2 H_2+\textrm{m}_3 P^\textrm{L}+\textrm{m}_4 \textrm{D}_1 \textrm{P}^\textrm{D}_1}{2 (H_1+H_2)} \\
  + \frac{\textrm{m}_5 \textrm{D}_2 \textrm{P}^\textrm{D}_2+\textrm{m}_6 R_1+\textrm{m}_7 R_2+\textrm{m}_8}{2 (H_1+H_2)} \leq \mbox{RoCoF}_{\textrm{max}} 
\end{multline}

Eq.~(\ref{eq:RoCoF_regional_regression}) is obtained by estimating term `$\textrm{A}_i \cdot \omega_i$' in (\ref{eq:RoCoF_analytical_TwoRegions}) using a linear combination of every system magnitude (i.e.~$H_1$, $H_2$, $P^\textrm{L}$, ...) divided by term `$2 (H_1+H_2)$', so as to obtain a linear formulation for the RoCoF constraint. The procedure for estimating the oscillation parameters $\textrm{A}_i$ and $\omega_i$ numerically (i.e.~finding the optimal value for the constant terms `$\textrm{m}$') is described in Algorithm~\ref{alg:MultiAreaSamples}, while an implementation of this algorithm using MATLAB and Simulink is freely available in \cite{GithubMultiArea}.

\begin{algorithm}[!t]
  \caption{Numerical estimation of terms `$\textrm{m}$' in eq.~(\ref{eq:RoCoF_regional_regression})}
 \label{alg:MultiAreaSamples}
 
 \begin{algorithmic}[1]
 
 \renewcommand{\algorithmicrequire}{\textbf{Input:}}
 \renewcommand{\algorithmicensure}{\textbf{Output:}}
 
 \REQUIRE range of operating conditions for the power system
 \ENSURE estimation of $\textrm{A}_i$ and $\omega_i$ for each condition 
 
  \FOR {several feasible values of $P^\textrm{L}$}
      \FOR {several splits of D among the regions}
      
          \STATE $H_\textrm{total} = P^\textrm{L}/(2 \cdot \textrm{RoCoF}_\textrm{max})$ \label{lineCode:CloseToBoundaryCOI_start}
          
          \STATE $R_\textrm{total} = H_\textrm{total}/k^*$  \label{lineCode:CloseToBoundaryCOI_k}
          
          \IF {$R_\textrm{total} < P^\textrm{L} - \textrm{D}\cdot\textrm{P}_\textrm{D}\cdot \Delta f_\textrm{max}$} \label{lineCode:qss}
              \STATE $R_\textrm{total} = P^\textrm{L} - \textrm{D}\cdot\textrm{P}_\textrm{D}\cdot \Delta f_\textrm{max}$ 
          \ENDIF \label{lineCode:CloseToBoundaryCOI_end}
          
          \FOR {several splits of $R_\textrm{total}$ among the regions} 
            \FOR {several splits of $H_\textrm{total}$ among the regions} 
                \STATE run dynamic simulation (e.g.~Simulink)
                
                \WHILE {RoCoF in the region $> \textrm{RoCoF}_\textrm{max}$} \label{lineCode:CloseToBoundaryRegion_start}
                    \STATE $H_\textrm{total}$ = $H_\textrm{total}$ + slight increase
                    \STATE run dynamic simulation
                \ENDWHILE \label{lineCode:CloseToBoundaryRegion_end}
                
                \STATE estimate $\textrm{A}_i$, $\omega_i$ (e.g.~\href{https://github.com/badber/TwoRegion_Frequency/blob/master/TwoRegion_RoCoF/Region1_faulted/estimate_oscillation.m}{\underline{`estimate\_oscillation.m'}})
                
                \STATE record features: system state ($H_1$, $H_2$, $P^\textrm{L}$, ...) 
                
                \STATE record labels: estimated `$\textrm{A}_i \cdot \omega_i$' 
                
            \ENDFOR 
          \ENDFOR 
      \ENDFOR
  \ENDFOR 
  
  \STATE regression with features \& labels (\href{https://github.com/badber/TwoRegion_Frequency/blob/master/TwoRegion_RoCoF/Region1_faulted/Rocof_regression.m}{\underline{`Rocof\_regression.m'}})
 
 \end{algorithmic} 
\end{algorithm}

\begin{figure}[!t]
	\centering
    \begin{subfigure}[b]{0.5\textwidth}
        \hspace*{-8mm}
        \centering
        \includegraphics[width=3.1in]{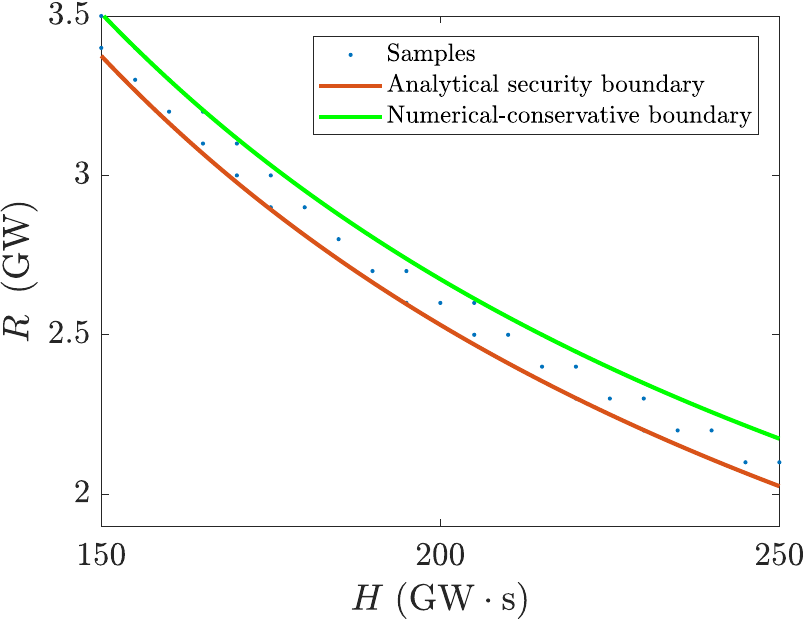}
        \captionsetup{width=.8\linewidth} 
        \caption{\footnotesize Using samples within 0.05Hz of the the nadir security boundary.}
        \vspace{3mm}
        \label{fig:SecurityBoundary_TooConservative}
    \end{subfigure}%
    ~ 
    \newline 
    \begin{subfigure}[b]{0.5\textwidth}
        \hspace*{-8mm}
        \centering
        \includegraphics[width=3.1in]{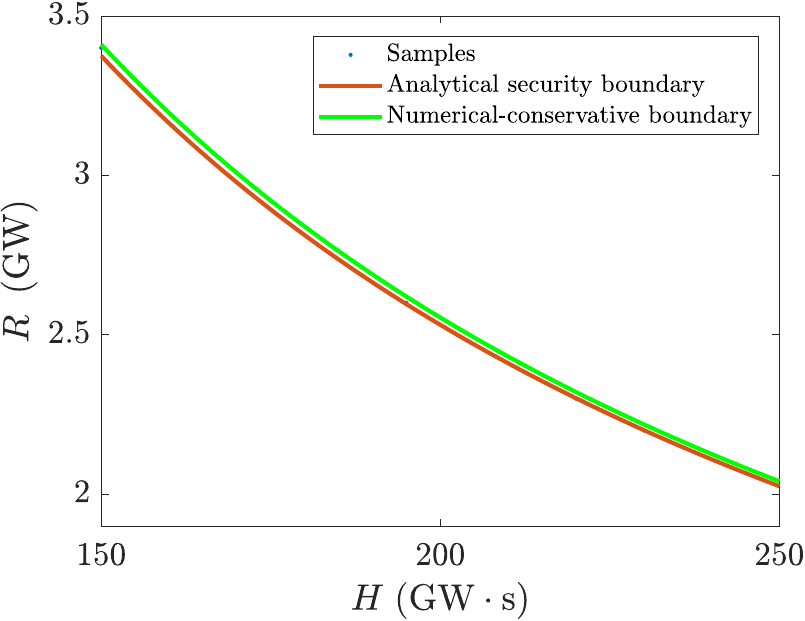}
        \captionsetup{width=.8\linewidth} 
        \caption{\footnotesize Using samples within 0.01Hz of the the nadir security boundary.}
        \label{fig:SecurityBoundary_good}
    \end{subfigure}
    \caption{Estimation of the COI frequency-stability boundary through a regression, using samples obtained from dynamic simulations.}
    \label{fig:SecurityBoundary_importanceSamples}
\end{figure}

Algorithm~\ref{alg:MultiAreaSamples} generates a number of feasible operating points for the power system, that exactly respect the regional RoCoF limit: dynamic simulations are run for a number of operating points, and samples that just barely respect the RoCoF requirement are recorded. 
Finally, a regression is run on these samples to find the right values for terms `$\textrm{m}$' in eq.~(\ref{eq:RoCoF_regional_regression}). Some further explanation on Algorithm~\ref{alg:MultiAreaSamples}:
\begin{itemize}
    \item Algorithm~\ref{alg:MultiAreaSamples} starts by considering a system operating point that just barely respects the COI frequency stability (i.e.~that neglects the inter-area oscillations), as described in lines \ref{lineCode:CloseToBoundaryCOI_start} through \ref{lineCode:CloseToBoundaryCOI_end}. Parameter `$k^*$' in line~\ref{lineCode:CloseToBoundaryCOI_k} refers to the condition for respecting the COI nadir requirement, deduced in \cite{FeiStochastic} as $H\cdot R\geq k^{*} $ (where $k^*$ is a function of the system operating state). 
    
    \item The reason for starting from a system operating point that just barely respects the COI frequency stability is to iteratively reach an operating point precisely at the regional RoCoF security boundary, using the loop in lines \ref{lineCode:CloseToBoundaryRegion_start} through \ref{lineCode:CloseToBoundaryRegion_end}. Using only samples at the RoCoF security boundary allows to obtain a more accurate regression for `$\textrm{A}_i \cdot \omega_i$', as explained later in this section and illustrated in Fig.~\ref{fig:SecurityBoundary_importanceSamples}.
    
    \item The several `for' loops included in the algorithm have the goal of generating samples that consider several possible splits of the system magnitudes ($H$, $R$, $\textrm{D}$...) among the regions. While there is no particular number of samples that is required, this algorithm allows to obtain a wide range of system operating points, appropriately covering the space of feasible operating points of the power system. 
    
\end{itemize}

The last line in Algorithm~\ref{alg:MultiAreaSamples} applies a regression to estimate the values of `$\textrm{A}_i \cdot \omega_i$', as this is the expression that needs to be estimated in eq.~(\ref{eq:RoCoF_analytical_TwoRegions}). The system operating points are used as features for the regression, which then will provide the values of terms `m' in eq.~(\ref{eq:RoCoF_regional_regression}). A conservative regression is used, so as to guarantee that the resulting constraints will lead to a frequency-stable operating point in every circumstance, with the tradeoff of a tighter frequency-security region if the security boundary is far from being linear. Nevertheless, as will be demonstrated in Section~\ref{Sec:CheckConservativeRoCoF_i}, this tradeoff is not significant for practical purposes in a power system.

The conservative regression can be computed solving the following constrained optimisation:
\begin{alignat}{3} 
& \underset{\theta}{\textrm{min}}  \quad \;\;
&& \frac{1}{2}\norm{\vphantom{\frac{\scriptstyle 1}{2}} X \cdot \theta - y}^2 \label{eq:ConstrainedLS} \\
& \; \text{s.t.} \quad \;
&& \quad X \cdot \theta \geq y \label{eq:ConstrainedLS2}
\end{alignat}
Where $y$ are the regression labels, $X$ are the regression features and $\theta$ is the vector of the regression parameters (in this case, $\theta = [\textrm{m}_1, ..., \textrm{m}_8]$ for eq.~\ref{eq:RoCoF_regional_regression}). The above optimisation problem simply defines a constrained least-squares regression, which forces the resulting regression to be above all training samples contained in vector $y$: in this case, the regression is forced to overestimate the regional RoCoF caused by the inter-area oscillation, i.e.~by term  `$\textrm{A}_i \cdot \omega_i$' in eqs.~(\ref{eq:RoCoF_analytical_TwoRegions}).

As explained before, Algorithm~\ref{alg:MultiAreaSamples} only uses samples close to the frequency-security boundary to train the conservative regression (i.e.~samples that barely respect the $\textrm{RoCoF}_\textrm{max}$ limit within the region). To illustrate the importance of training the conservative regression using only samples that fall very close to the frequency-security boundary, consider Fig.~\ref{fig:SecurityBoundary_importanceSamples}: as a conservative above-all-samples regression is used, training it only with samples close to the security boundary avoids an overly-conservative RoCoF constraint. Note that this Fig.~\ref{fig:SecurityBoundary_importanceSamples} is for illustration purposes only, as it has been generated by considering the COI nadir limit ($H\cdot R\geq k^{*} $, as deduced in \cite{FeiStochastic}): since the analytical expression for that limit can be obtained, it is possible to visually understand the advantage of using appropriate samples when estimating the boundary through a regression.

The splits of the system parameters considered in Algorithm~\ref{alg:MultiAreaSamples} must be within the typical operating ranges of the power system, as the goal is to implement the resulting constraints in a scheduling algorithm that will consider precisely these ranges. Note that this algorithm is only run once and this is done `offline’, i.e.~before the actual system scheduling is run: the system operator would only have to re-run the algorithm when the generation mix changes significantly, possibly on an annual basis. The algorithm should be re-run once renewable penetration has increased sufficiently (which would shift the ranges of inertia in the system noticeably), as the regressions have to be trained for the ranges considered in the Unit Commitment.

Finally, note that the model for Primary Frequency Response used in the dynamic simulations is the ramp described in eq.~(10) of Part I. This is due to the fact that the focus of this work is to obtain frequency-stability constraints for the scheduling problem, for which system operators define the requirements for PFR to be fully delivered by a pre-determined time. In the case of Great Britain, the requirement is to provide full response by 10s after a fault: all machines participating in the PFR service will adjust their droop gain in order to comply with this speed. Our previous work \cite{LuisMultiFR} showed how the ramp for PFR represents accurately the droop control for frequency response from synchronous generators, while providing simple market rules that all providers of PFR can comply with.

\subsection{Conservativeness of the resulting RoCoF constraint} \label{Sec:CheckConservativeRoCoF_i}

\begin{figure} [!t]
\hspace*{-8mm}
    \centering
    \includegraphics[width=3.1in]{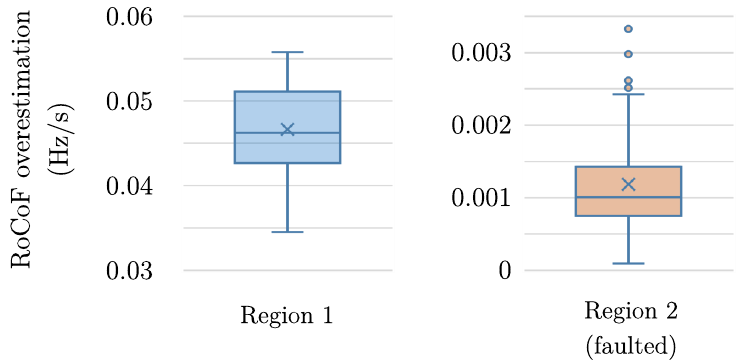}
    \vspace{1mm}
    \caption{Assessment of conservativeness in regional RoCoF constraints.}
    \label{fig:ConservativenessRocof}
\end{figure}

In this section we demonstrate that the resulting RoCoF constraints do indeed respect the $\textrm{RoCoF}_\textrm{max}$ limit in each region. In order to assess how the resulting RoCoF constraints perform, a number of system operating states was generated that exactly meet the RoCoF constraint in each region, i.e.~that exactly meet constraint~(\ref{eq:RoCoF_regional_regression}). Then, each system operating point was fed into a dynamic simulation to compute the actual RoCoF that would occur for that operating point. The comparison of the computed RoCoF and the predicted RoCoF using constraint~(\ref{eq:RoCoF_regional_regression}) in each region is presented in the box-and-whiskers plots in Fig.~\ref{fig:ConservativenessRocof}.

The results in Fig.~\ref{fig:ConservativenessRocof} were obtained from 270 samples of system operating points in the ranges $H \in [110,170]\textrm{GW}\cdot \textrm{s}$, $R \in [1.6,3]\textrm{GW}$, $P^{\textrm{L}} \in [1.2,1.8]\textrm{GW}$ and $\textrm{P}_\textrm{D} \in [20,45]\textrm{GW}$, ranges within the GB 2030 system discussed in Section~\ref{sec:GB2030}. The RoCoF limit in both regions was set to $\textrm{RoCoF}_\textrm{max}=1\textrm{Hz/s}$ and the damping factor to $\textrm{D}=0.5\%/\textrm{Hz}$, while the transmission corridor considered had a line reactance of $\textrm{X}_{1,2}=50\Omega$ and voltages of $\textrm{V}_{1}=\textrm{V}_{2}=400\textrm{kV}$. 
The box-and-whiskers plot represents 75\% of the distribution within the box, while the whole distribution is contained within the whiskers (excluding any outliers). The dots outside the whiskers represent outliers (i.e.~samples placed 1.5 times the inter-quantile range outside the limits of the box). The median is represented by the horizontal line inside the model, while the `x' represents the mean.

As shown in Fig.~\ref{fig:ConservativenessRocof}, the proposed RoCoF constraints are almost perfect for region 2, while they overestimate RoCoF by 0.055Hz/s in the worst case for region 1 (with a median overestimation of 0.046Hz/s). These results are consistent with the analytical deduction of the RoCoF constraint: 
the proposed constraint neglects the effect of the attenuation of inter-area oscillations, which entails no approximation for the faulted region but a certain approximation for other regions. 
The maximum RoCoF in the faulted region occurs at the very instant of the fault, when there is no attenuation of oscillations present; however, the maximum RoCoF in the non-faulted region 1 occurs a few instants afterwards, so neglecting the attenuation entails some conservativeness.

\subsection{Nadir constraint: numerical estimation} \label{sec:NadirRegionalNumericalEstimation}

In order to obtain a linear formulation for the nadir constraints, a numerical estimation must be applied to the integrals of post-fault frequency deviation appearing in them. In a similar way as done in Section~\ref{sec:RoCoFnumericalEstimation} for the RoCoF constraint, here we use a linear combination of every system magnitude to estimate these integrals:

\begin{multline}
\int_{0}^{\textrm{t}_\textrm{nadir}}\int_{0}^{t}\left[\Delta f_1(\tau)-\Delta f_2(\tau)\right] \textrm{d}\tau \ \textrm{d}t \quad \xrightarrow{\quad \textrm{estimated by} \quad } \\ 
f(\textrm{system operating state})=  \textrm{m}'_1 H_1+\textrm{m}'_2 H_2+\textrm{m}'_3 P^\textrm{L} \\
+\textrm{m}'_4 \textrm{D}_1 \textrm{P}^\textrm{D}_1+\textrm{m}'_5 \textrm{D}_2 \textrm{P}^\textrm{D}_2+\textrm{m}'_6 R_1+\textrm{m}'_7 R_2+\textrm{m}'_8
\end{multline}

\begin{multline}
\int_{0}^{\textrm{t}_\textrm{nadir}}\Delta f_i(t) \ \textrm{d}t \quad \xrightarrow{\quad \textrm{estimated by} \quad } \\
f(\textrm{system operating state})=  \textrm{m}''_1 H_1+\textrm{m}''_2 H_2+\textrm{m}''_3 P^\textrm{L} \\
+\textrm{m}''_4 \textrm{D}_1 \textrm{P}^\textrm{D}_1+\textrm{m}''_5 \textrm{D}_2 \textrm{P}^\textrm{D}_2+\textrm{m}''_6 R_1+\textrm{m}''_7 R_2+\textrm{m}''_8
\end{multline}

The procedure for performing this numerical estimations is equivalent to Algorithm~\ref{alg:MultiAreaSamples}, and an implementation of this algorithm is also available in the code repository \cite{GithubMultiArea}. In the following section, the conservativeness introduced by these numerical estimations is quantified, since the constrained least-squares regression defined in (\ref{eq:ConstrainedLS}) and (\ref{eq:ConstrainedLS2}) is also used here. Using this constrained regression guarantees that the frequency nadir will be above the stability threshold of $\Delta f_\textrm{max}$ in every circumstance.

\subsection{Conservativeness of the resulting nadir constraint} \label{sec:NadirRegionalAssessment}

\begin{figure} [!t]
\hspace*{-7mm}
    \centering
    \includegraphics[width=2.8in]{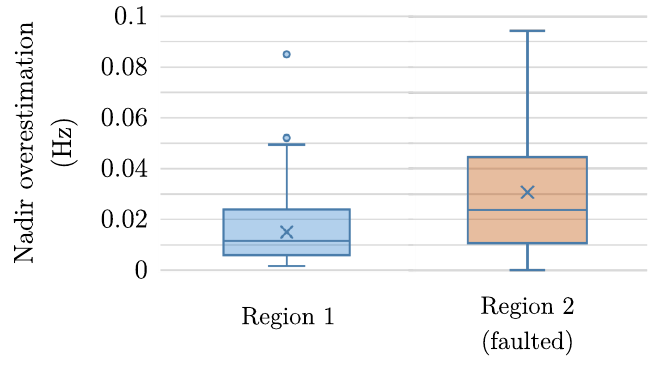}
    \caption{Assessment of conservativeness in regional nadir constraints.}
    \label{fig:ConservativenessNadir}
\end{figure}

In order to assess how the resulting nadir constraints for each region perform, a number of system operating states was generated that exactly meet the nadir constraint in each region. 
A dynamic simulation was then run for each of these states, so as to compute the actual nadir that would occur for that operating point. Fig.~\ref{fig:ConservativenessNadir} presents the comparison of the computed nadirs from the dynamic simulations and the predicted nadirs using the proposed constraint in each region.

The results in Fig.~\ref{fig:ConservativenessNadir} were obtained from 300 samples of system operating points in the ranges $H \in [60,125]\textrm{GW}\cdot \textrm{s}$, $R \in [3,5.2]\textrm{GW}$, $P^{\textrm{L}} \in [1.2,1.8]\textrm{GW}$ and $\textrm{P}_\textrm{D} \in [20,45]\textrm{GW}$. The nadir limit in both regions was set to $\Delta f_\textrm{max}=0.8\textrm{Hz}$ and the damping factor to $\textrm{D}=0.5\%/\textrm{Hz}$, while the transmission corridor considered had a line reactance of $\textrm{X}_{1,2}=50\Omega$ and voltages of $\textrm{V}_{1}=\textrm{V}_{2}=400\textrm{kV}$. The results show that nadir is on average overestimated by under 0.02Hz in region 1 (i.e.~a 2\% conservativeness) and 0.03Hz in region 2 (3.75\% conservativeness), while the maximum overestimation is of 0.09Hz for either region.

\section{Case studies} \label{Sec:CaseStudies}

The simultaneous optimisation of energy and ancillary services for frequency support can be achieved by appropriately constraining a scheduling algorithm to guarantee frequency stability. In this section we demonstrate, through several case studies, the implications of distinct regional frequencies for the procurement of ancillary services.  

It could be argued that inter-area oscillations are typically managed through oscillation dampers such as Power Systems Stabilizers (PSS), and therefore the scheduling algorithm need not be constrained to consider these oscillations. However, the contribution of these devices will be very limited, if not completely negligible, for time-scale of interest for the RoCoF and nadir after a generation loss: the highest RoCoF occurs in the very instant of the outage, and therefore a PSS cannot reduce the regional RoCoF. In the case of the frequency nadir, this can occur within just a few seconds after a fault in power grids with low inertia, when the contribution of a PSS will be limited. Therefore, the Unit Commitment must be constrained to guarantee post-fault frequency security in systems exhibiting inter-area oscillations, as sufficient inertia and frequency response must be available to contain them.

\subsection{Frequency-secured stochastic scheduling model} \label{sectionSUC}

The tool used here for conducting simulations of a power system scheduling is a two-stage Stochastic Unit Commitment (SUC) model. The Unit Commitment is a problem of particular interest when studying the provision of frequency services, since the commitment state of synchronous generators determines the level of system inertia, the key driver of post-fault frequency evolution.

This SUC solves an optimisation problem in which the decision variables are the commitment status (on/off) of each generator and the power output of online generators. The objective is to minimise the system's expected operational cost, subject to several constraints: generation-demand balance each hour, generation limits for each unit, inter-temporal constraints such as start-up times or maximum ramp rates of thermal units, etc. (for a full description of these constraints, refer to \cite{LuisCovid}). The SUC is solved with a 24h lookahead, by building a forecast for demand and RES generation in the following 24h and computing the schedule for dispatchable generators with an hourly resolution.

\begin{figure}[!t]
    \centering
    \includegraphics[width=3.1in]{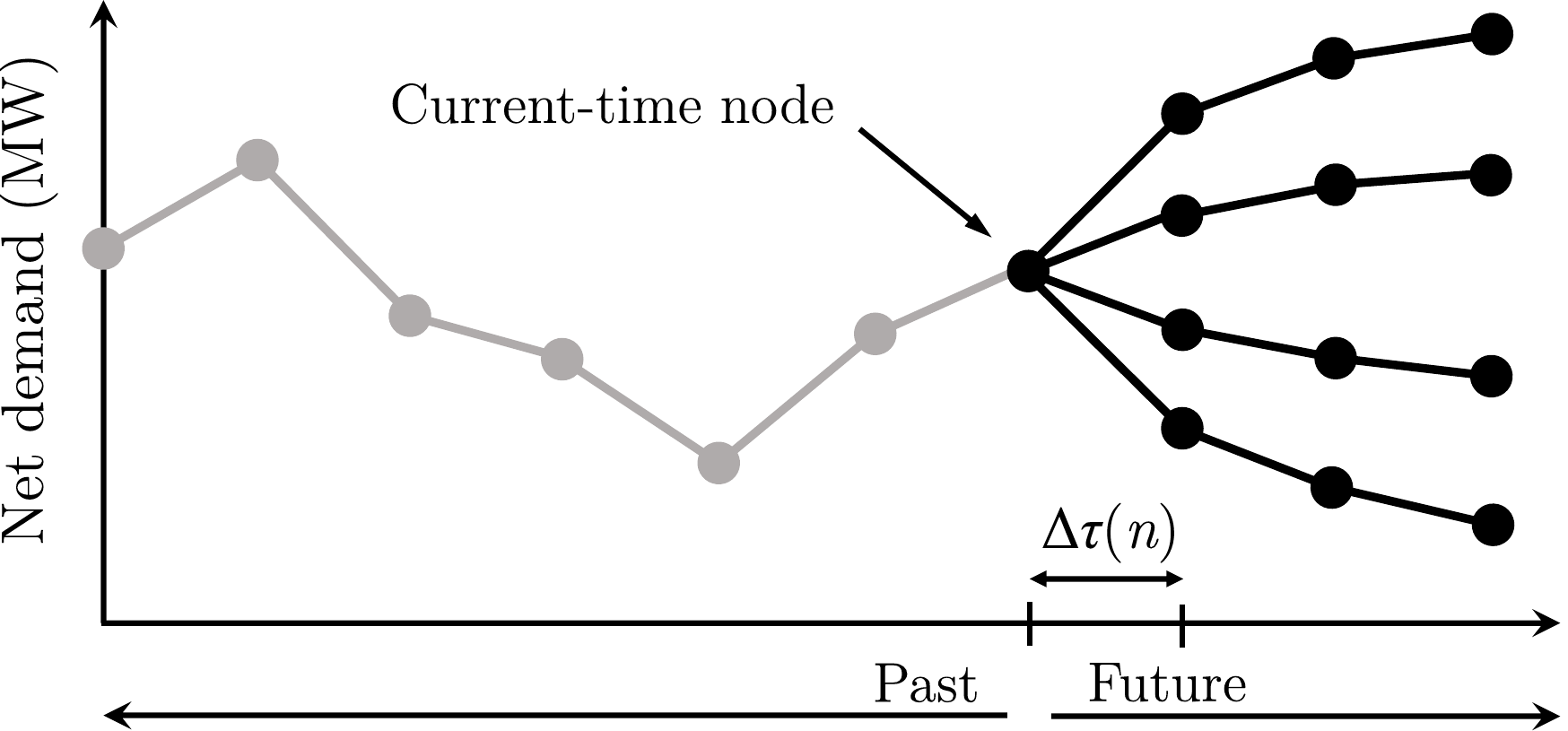}
    \vspace{3mm}
    \caption{Scenario tree used in the Stochastic Unit Commitment \cite{LuisEFR}. 
    }
    \label{fig:ScenarioTree}
\end{figure}

In order to account for the uncertainty in renewable generation, the SUC considers a number of possible realisations of the stochastic variable in the system, i.e.~RES generation. The infinite set of possible realisations is discretised into a representative scenario tree such as the one in Fig.~\ref{fig:ScenarioTree}, 
in which each node $n$ has a certain probability of occurrence.
The objective function of the SUC, which is formulated as a Mixed-Integer Linear Program (MILP), is the expected operation cost over all nodes in the scenario tree:
\begin{equation} \label{objectiveSUC}
\text{min}\quad \sum_{n \in \mathcal{N}}\pi(n) \left( \sum_{g \in \mathcal{G}}C_g(n) + \textrm{c}^\textrm{LS}\cdot \Delta\tau(n) \cdot P^\textrm{LS}(n) \right)
\end{equation}
The sum of the operating cost $C_g$ of all thermal units at each node $n$, plus the penalties paid for load shedding, is weighted by the probability of reaching that node, $\pi(n)$. 
The operating cost of a group of generating units $g$ is given by:
\begin{equation} 
C_g(n)=\textrm{c}_g^\textrm{st}\cdot N_g^\textrm{sg}(n)+\Delta\tau(n)\left[\textrm{c}_g^\textrm{nl}\cdot N_g^\textrm{up}(n)+\textrm{c}_g^\textrm{m}\cdot P_g(n)\right]
\end{equation}
Note that generating units with the same characteristics are clustered in the SUC to reduce the computational burden, as proposed in \cite{AlexEfficient}.

The regional frequency constraints obtained in Section~\ref{Sec:Applicability} are also included in the SUC, which implies that the scheduling solution is guaranteed to contain enough resources (inertia, PFR) so that the real-time controls can maintain the frequency stability of the system after a contingency. This implies that the SUC co-optimises energy and ancillary services: if the inertia provided by generators needed online to produce energy is not sufficient to comply 
with the frequency constraints, the SUC will commit a higher number of generators. The same reasoning applies to PFR, which is provided by the headroom in thermal units, as defined by:
\begin{equation}
    R_g \leq \textrm{R}_g^\textrm{max}
\end{equation}
\begin{equation}
    R_g \leq \textrm{R}_g^\textrm{slope}\cdot \left( N_g^\textrm{up}\cdot\textrm{P}_g^\textrm{max} - P_g \right)
\end{equation}
The total system inertia is the sum of inertia provided by all committed generators, minus the inertia lost from the outaged generator: 
\begin{equation} \label{eq:DefinitionSystemH}
H=\sum_{g \in \mathcal{G}}\textrm{H}_g\cdot \textrm{P}_g^{\textrm{max}}\cdot N_g^{\textrm{up}} -  \textrm{H}_\textrm{L} \cdot \textrm{P}^{\textrm{max}}_{\textrm{L}} 
\end{equation}

The cost of these ancillary services is implicit in the cost of energy, as part-loaded generators operate at a lower efficiency, which increases their per-MWh-cost. In other words, the same amount of energy produced by a small number of fully-loaded generators is cheaper than using a higher number of part-loaded generators, while the latter case increases the inertia and PFR in the system.

Note that there is an additional implicit cost associated to the provision of ancillary services during low net-demand conditions: when RES output is high, covering most if not all of the demand at a given time, a number of thermal units must be kept online to provide inertia and PFR. As these plants cannot operate below their minimum stable generation point $\textrm{P}_{\textrm{g}}^{\textrm{msg}}$, the energy generated by these plants displaces energy that could be generated from RES. This RES curtailment is purely due to guaranteeing system stability in the event of a contingency, and therefore the increase in energy costs is effectively the cost of ancillary services during these periods.

In the simulation results presented in this Section~\ref{Sec:CaseStudies}, the SUC is run for one year of operation of the system. A scenario tree with 24h lookahead is built in every hour of the year. Only the scheduling decision for the current-time node is recorded, with all other decisions for the future 24h being discarded before building the new scenario tree. We refer to this approach as `rolling planning', which is similar to the `model predictive control' strategy used in the control theory community. By simulating the operation of the power system during a full year, it is possible to understand the need for ancillary services required to guarantee frequency stability while accounting for the different RES and demand levels throughout a year.

\subsection{Great Britain 2030 power system} \label{sec:GB2030}

The power system considered to conduct simulations in this paper is a feasible representation of the Great Britain system by 2030, which considers a partially-decarbonised generation mix characterised by: 1) a nuclear fleet providing baseload and driving the \textit{N}-1 reliability requirement for frequency services, due to the large power rating of some of these generators; 2) a fleet of gas-fired power plants including Combined-Cycle Gas Turbines (CCGTs) and Open-Cycle Gas Turbines (OCGTs); 
and 3) a high increase in renewable capacity from current levels, with 60GW of wind power considered in this paper (corresponding to the average wind capacity of the four scenarios envisioned for 2030 by the British system operator \cite{NationalGridFES2020}). 
The parameters defining the operational characteristics of the thermal fleet are included in Table~\ref{TableThermal_intro}. 
    
Electricity demand ranges from 20GW to a peak of 60GW, accounting for daily and seasonal variations,
with 90\% of the load located in England and 10\% in Scotland \cite{DemandEnglandScotland}. Wind capacity is split by 50\%-50\% between England and Scotland \cite{BEISwind2019}, and the transmission corridor connecting these two regions has a voltage of 400kV, line reactance of $50\Omega$ and thermal limit of 7.5GW \cite{ETYS2020}. 
The England-Scotland regions were identified by \cite{InmaMultiArea,GBtwoRegions} as an appropriate division of the British grid into areas where synchronous machines show a similar frequency behaviour.

\begin{table}[!t]
\captionsetup{justification=centering, textfont={sc,footnotesize}, labelfont=footnotesize, labelsep=newline} 
\renewcommand{\arraystretch}{1.6}
\caption{Characteristics of thermal plants in GB's 2030 system}
\label{TableThermal_intro}
\centering
\begin{tabular}{m{3.1cm}|m{0.65cm}m{0.6cm}m{0.7cm}|m{0.6cm}m{0.65cm}}
    \multicolumn{1}{l}{} & \multicolumn{3}{c|}{England} & \multicolumn{2}{c}{Scotland}\\ \cline{2-6}
    
    \multicolumn{1}{l}{} & Nuclear & CCGT & OCGT & CCGT & OCGT\\ 
\hline
\hspace*{-2mm} \raggedright Number of Units & 4 & 80 & 25 & 20 & 5 \\ \hline
\hspace*{-2mm} \raggedright Rated Power $\textrm{P}_{\textrm{g}}^{\textrm{max}}$ (MW) & 1800 & 500 & 100 & 500 & 100\\ \hline
\hspace*{-2mm} \raggedright Min Stable Gen $\textrm{P}_{\textrm{g}}^{\textrm{msg}}$ (MW) & 1800 & 250 & 50 & 250 & 50\\ \hline
\hspace*{-2mm} \raggedright No-Load Cost $\textrm{c}^{\textrm{nl}}_g$ (\pounds/h) & 0 & 4500 & 3000 & 4500 & 3000\\ \hline
\hspace*{-2mm} \raggedright Marginal Cost $\textrm{c}^{\textrm{m}}_g$ (\pounds/MWh) & 10 & 46 & 200 & 46 & 200\\ \hline
\hspace*{-2mm} \raggedright Start-up Cost $\textrm{c}^{\textrm{st}}_g$ (\pounds) & N/A & 10000 & 0 & 10000 & 0\\ \hline
\hspace*{-2mm} \raggedright Start-up Time (h) & N/A & 4 & 0 & 4 & 0\\ \hline
\hspace*{-2mm} \raggedright Min Up Time (h) & N/A & 4 & 0 & 4 & 0\\ \hline
\hspace*{-2mm} \raggedright Min Down Time (h) & N/A & 1 & 0 & 1 & 0\\ \hline
\hspace*{-2mm} \raggedright Inertia Constant $\textrm{H}_g$ (s) & 5 & 5 & 5 & 5 & 5\\ \hline
\hspace*{-2mm} \raggedright PFR capacity $\textrm{R}_g^\textrm{max}$ (MW) & 0 & 50 & 20 & 50 & 20\\ \hline
\hspace*{-2mm} \raggedright PFR slope $\textrm{R}_g^\textrm{slope}$ & 0 & 0.5 & 0.5 & 0.5 & 0.5\\ \hline
\hspace*{-2mm} \raggedright Emissions (g $\mbox{CO}_{2}$/kWh) & 0 & 368 & 833 & 368 & 833\\
\hline
\end{tabular}
\end{table}

Regarding post-fault frequency limits, the values defined in GB are \cite{StandardNationalGrid}: 
\begin{itemize}
    \item RoCoF must be below 0.125Hz/s at all times to avoid the tripping of RoCoF-sensitive protection relays (i.e.~$\textrm{RoCoF}_\textrm{max}=0.125\textrm{Hz/s}$).
    \item The frequency nadir must never be below 49.2Hz to prevent the activation of Under Frequency Load Shedding (i.e.~$\Delta f_\textrm{max}=0.8\textrm{Hz}$).
    \item Frequency must recover to be above 49.5Hz within 60s after the outage, referred to as the frequency quasi-steady-state requirement (i.e.~$\Delta f_\textrm{max}^\textrm{ss}=0.5\textrm{Hz}$).
\end{itemize}
Due to increasing penetration of non-synchronous RES in GB, which would significantly increase the procurement cost of frequency services needed to comply with this regulation, the RoCoF limit is in the process of being relaxed to 1Hz/s in the whole network \cite{LossOfMainsProgramme}, therefore this paper considers a value of $\textrm{RoCoF}_\textrm{max}=1\textrm{Hz/s}$.

Since the problem solved here is the Unit Commitment, which focuses on obtaining the commitment decision for the different thermal generators, the only network constraints considered are the thermal limit for transfer of power between regions, and the regional frequency stability constraints. In practice, an Optimal Power Flow could solved be once the commitment decision of the generators is fixed, which considers a full representation of the electricity network.

\subsection{Uniform-frequency model vs. Regional-frequency model} \label{Sec:SingleAreaVsMultiArea}

\begin{figure} [!t]
\hspace*{-2mm}
    \centering
    \includegraphics[width=3.3in]{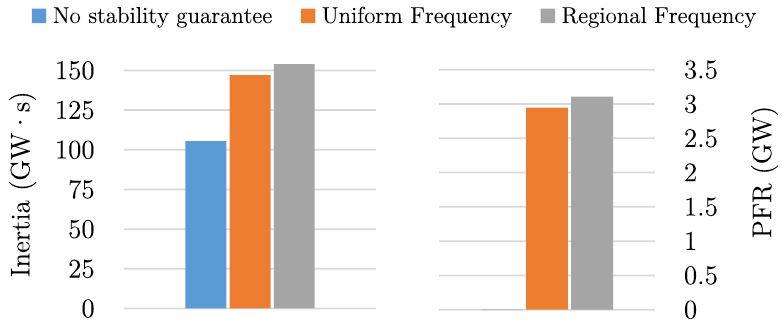}
    \vspace{1mm}
    \caption{Average hourly volumes of ancillary services needed to comply with post-fault frequency limits.}
    \vspace{1mm}
    \label{fig:SingleAreaVsMulti}
\end{figure}

\begin{figure} [!t]
\hspace*{-5mm}
    \centering
    \includegraphics[width=3.3in]{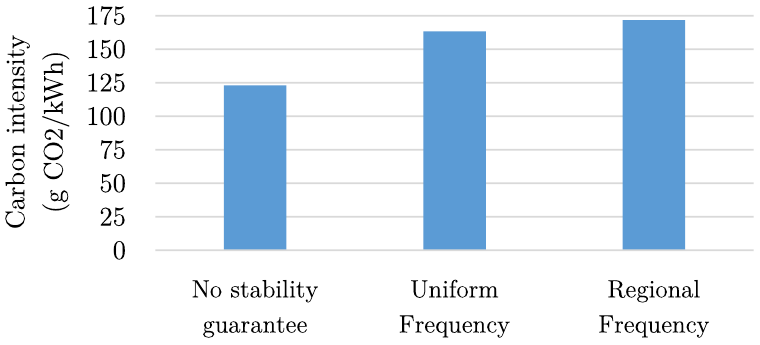}
    \vspace{1mm}
    \caption{Carbon intensity resulting from the optimal scheduling of dispatchable generation.}
    \label{fig:CarbonIntensity}
\end{figure}

We first analyse the consequences of guaranteeing frequency stability stability in the system scheduling, and study the extra need for ancillary services in the presence of inter-area oscillations in frequency. The SUC is secured against a generation loss of 1.8GW occurring in England, i.e.~$P^\textrm{L}=1.8\textrm{GW}$.

Fig.~\ref{fig:SingleAreaVsMulti} compares three cases: 1) an energy-only scheduling, obtained from running the SUC with no frequency constraints implemented (therefore the scheduling solution does not guarantee frequency stability in the event of an outage); 2) a frequency-secured SUC, considering the uniform frequency model proposed in \cite{FeiStochastic}; and 3) a frequency-secured SUC, considering the constraints for regional frequency stability proposed in this paper. In the first case, referred to as `No stability guarantee' in Fig.~\ref{fig:SingleAreaVsMulti}, the inertia present is simply a by-product of energy, provided by nuclear plants and gas-fired plants used during periods of low wind generation. The volume of PFR procured in this case is 0MW, since PFR is never a by-product of energy, as running part-loaded thermal plants would unnecessarily increase the cost of energy. In practice, most system operators mandate thermal units to provide PFR when generating, if their nominal capacity is larger than a predefined value; however, this scheduling with `No stability guarantee' does not consider such grid-code constraints, as it simply minimises the cost of producing energy.

On the other hand, both inertia and PFR volumes increase for cases that guarantee frequency stability. The results show that accounting for the regional variations in frequency would not significantly increase these volumes as compared to the `Uniform frequency' case, which neglects the inter-area oscillations in frequency: the inertia procured is around $7\textrm{GW}\cdot \textrm{s}$ higher (i.e.~5\% higher), while PFR increases by 150MW (that is, just above 5\% higher). In following sections of this paper, we however demonstrate that the location of inertia and PFR in a particular region is key to guarantee stability.

It is also insightful to visualize the carbon intensity of the system, presented in Fig.~\ref{fig:CarbonIntensity} for all three cases. The results show an increase of 40g$\mbox{CO}_{2}$/kWh from case `No stability guarantee' compared to case `Uniform frequency': since thermal plants must be committed to provide inertia and response, this increase in emissions is associated purely with stability actions that would have to be taken by the system operator. An additional 8.5g$\mbox{CO}_{2}$/kWh would be needed for guaranteeing stability in the presence of inter-area oscillations between England and Scotland.

\subsection{Where to procure inertia and response?} \label{Sec:LocationInertiaResponse}

\begin{figure} [!t]
\hspace*{-5mm}
    \centering
    \includegraphics[width=3.3in]{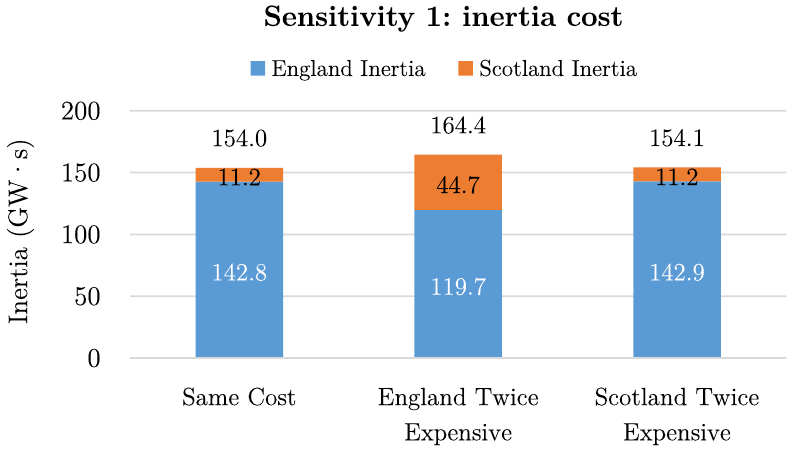}
    \vspace*{1mm}
    \caption{Average inertia procured in each region, for a 1.8GW loss occurring in England and a sensitivity analysis for the cost of inertia in each region.}
    \vspace*{1mm}
    \label{fig:InertiaCost_InertiaVolumes}
\end{figure}

\begin{figure} [!t]
\hspace*{-5mm}
    \centering
    \includegraphics[width=3.3in]{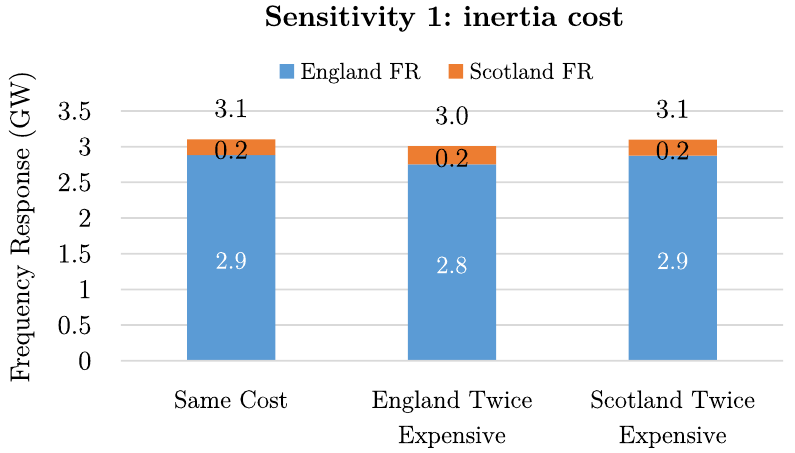}
    \vspace*{1mm}
    \caption{Average PFR procured in each region, for a 1.8GW loss occurring in England and a sensitivity analysis for the cost of inertia in each region.}
    \vspace*{1mm}
    \label{fig:InertiaCost_FRVolumes}
\end{figure}

This section has the aim of answering the following question: depending on the cost of inertia and response in each region, \textit{where should these services be procured?} To do so, we run a sensitivity analysis by adding an explicit cost penalty to inertia and response in each region,
which in the base case is of $\textrm{\pounds}250/\textrm{MW}$ for PFR and $\textrm{\pounds}5/\textrm{MW}\cdot \textrm{s}$ for inertia.
These values roughly correspond to the cost of procuring these services during periods of low net-demand, when the need for frequency services causes RES curtailment, as explained in Section~\ref{sectionSUC}. Then, this base case penalty is doubled in one of regions to understand how the SUC solution changes. Again, the SUC is secured against a generation loss of 1.8GW occurring in England.

The first sensitivity considered is the penalty for inertia in each region, with the results presented in Figs.~\ref{fig:InertiaCost_InertiaVolumes} and \ref{fig:InertiaCost_FRVolumes}. These results show that inertia is mostly located in England, since the 1.8GW fault takes place there and therefore inertia located in that region is the most effective means to contain RoCoF in England.
Fig.~\ref{fig:InertiaCost_InertiaVolumes} also illustrates that, when inertia is twice as expensive in England, the volume of inertia procured in England decreases, and this decrease is compensated by a higher volume of inertia procured in Scotland. Note that the total system inertia increases, since inertia in Scotland is less effective to contain the loss in England and therefore more inertia overall is needed. Fig.~\ref{fig:InertiaCost_FRVolumes} shows no significant changes on the volumes of PFR procured.

\begin{figure} [!t]
\hspace*{-5mm}
    \vspace*{1mm}
    \centering
    \includegraphics[width=3.3in]{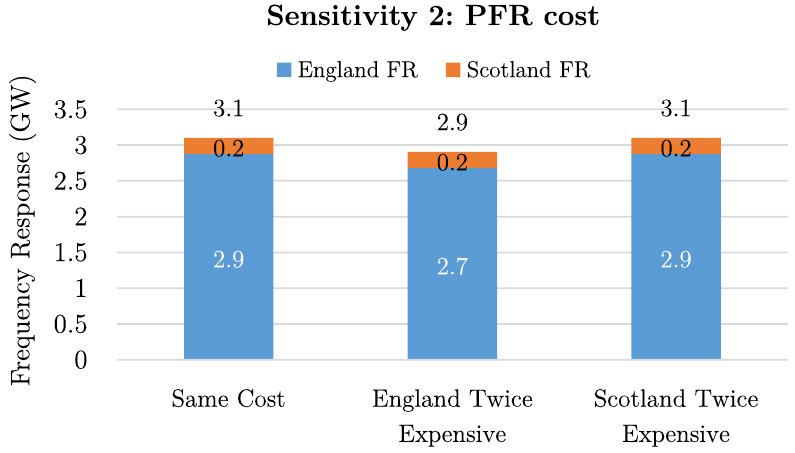}
    \caption{Average PFR procured in each region, for a 1.8GW loss occurring in England and a sensitivity analysis for the cost of PFR in each region.}
    \vspace*{1mm}
    \label{fig:FRCost_FRVolumes}
\end{figure}

\begin{figure} [!t]
\hspace*{-5mm}
    \centering
    \includegraphics[width=3.3in]{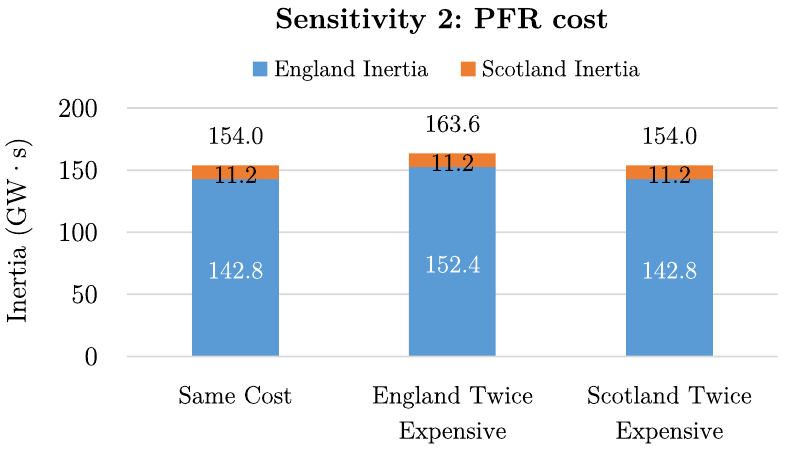}
    \caption{Average inertia procured in each region, for a 1.8GW loss occurring in England and a sensitivity analysis for the cost of PFR in each region.}
    \label{fig:FRCost_InertiaVolumes}
\end{figure}

The second sensitivity is the penalty for PFR in each region, presented in Figs.~\ref{fig:FRCost_FRVolumes} and \ref{fig:FRCost_InertiaVolumes}. Again, most of the inertia is procured in England, while Fig.~\ref{fig:FRCost_FRVolumes} shows that, if PFR is twice as expensive in England, a lower volume of response is procured in this region (2.7GW compared to 2.9GW). In turn, to compensate this drop in PFR and be able to contain the frequency drop in England, the volume of inertia in England increases ($152\textrm{GW}\cdot \textrm{s}$ compared to $142\textrm{GW}\cdot \textrm{s}$). The PFR procured in Scotland does not increase even though it is in this case cheaper than in England, which reflects the fact that PFR from Scotland must travel through the transmission corridor to contain the frequency nadir in England, and therefore it is less effective for this case of a loss occurring in England. 

The results in Figs.~\ref{fig:InertiaCost_InertiaVolumes} through \ref{fig:FRCost_InertiaVolumes} show that the volumes of inertia and PFR procured in Scotland simply do not drop below the minimum values to contain the frequency drop in that region. The location of the fault has a clear impact on where the frequency services are needed, as the results presented in this Section have shown that most inertia and response must be located in England if a large generation outage occurs in that region. However, the regional cost of frequency services also has a noticeable impact on the optimal volumes procured in each region. In the next section, we further analyse the impact of the fault location by considering a fault in Scotland, the low-inertia region given its excess wind generation.

\subsection{Impact of fault location: fault in the low-inertia region} \label{Sec:FaultScotland}

Here we analyse the implications of a generation loss taking place in the low-inertia region. We consider a 0.8GW loss in Scotland, corresponding to half of the capacity of the double-circuit HVDC interconnector named North Sea Link, expected to be commissioned in 2021 \cite{NorthSeaLink}.
Three cases are considered: 1) the uniform frequency model for frequency stability from \cite{FeiStochastic}, with a thermal limit between England and Scotland of 7.5GW; 2) the model for guaranteeing regional frequency stability proposed in this paper, with the same thermal limit; and 3) the model for regional frequency stability, but removing the thermal limit to understand the implications in has on procurement of frequency services in each region.

\begin{figure} [!t]
\hspace*{-5mm}
    \centering
    \includegraphics[width=3.3in]{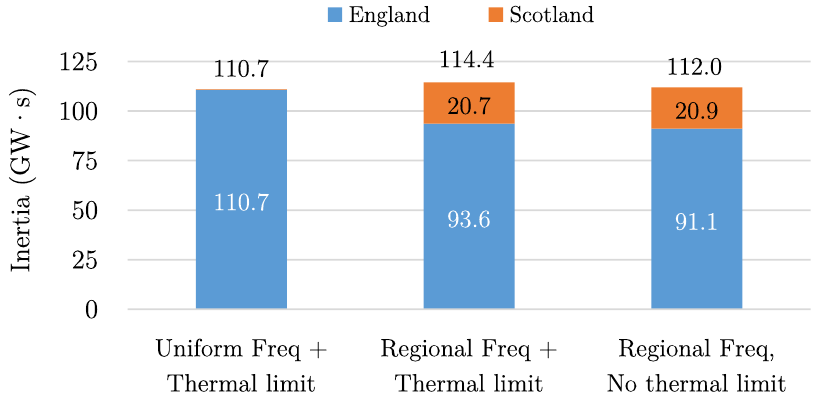}
    \vspace{1mm}
    \caption{Average inertia procured in each region throughout the year, for a 0.8GW fault taking place in Scotland.}
    \vspace{1mm}
    \label{fig:FaultScotlandInertia}
\end{figure}

The average inertia scheduled by the SUC in each region is included in Fig.~\ref{fig:FaultScotlandInertia}, for each of these three cases. By comparing the uniform-frequency model solution with the solution considering regional frequency stability, it is clear that some thermal plants must be committed in Scotland to provide inertia, as otherwise even this medium-size fault in that region would lead to violation of the RoCoF limit: the uniform-frequency model schedules no inertia in Scotland as it assumes that the location of inertia is irrelevant, given that frequency is considered a system-wide magnitude; however, zero inertia in Scotland would make frequency in that region drop virtually infinitely at the time of the fault, therefore violating the RoCoF limit.

By comparing the solutions respecting regional frequency stability with and without enforcing the thermal limit in the transmission corridors, it is demonstrated that this limit has the consequence of slightly increasing overall inertia in the system: since the number of thermal plants generating in Scotland is constrained by the thermal limit, a higher volume of inertia is procured in England. But since the loss occurs in Scotland, inertia procured in England is less effective to contain the frequency drop, and therefore a higher volume must be procured.

\begin{figure} [!t]
\hspace*{-5mm}
    \centering
    \includegraphics[width=3.3in]{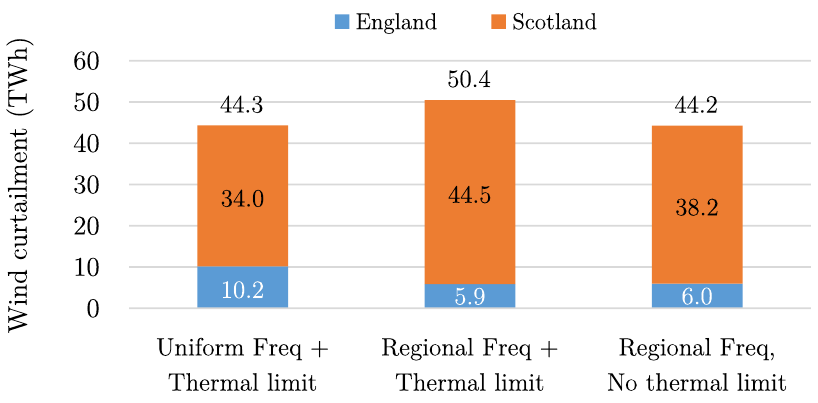}
    \vspace{1mm}
    \caption{Annual wind curtailment in each region, for a 0.8GW fault taking place in Scotland.}
    \label{fig:FaultScotlandWindCurtail}
\end{figure}

Furthermore, Fig.~\ref{fig:FaultScotlandWindCurtail} presents the impact on wind curtailment caused by the stability actions. Given the thermal limit of 7.5GW for the transmission corridors between Scotland and England, and that only 10\% of the total GB demand is consumed in Scotland, there is a maximum volume of energy that can be generated in Scotland so that the energy-export limit is respected. This has the effect of causing wind curtailment in Scotland, even in the `uniform frequency' case where no thermal generation is committed in Scotland. However, respecting regional frequency stability in Scotland aggravates this effect, as the thermal plants that must be online in that region to provide inertia displace wind generation so that the thermal limit is respected. The last column in Fig.~\ref{fig:FaultScotlandWindCurtail} shows that, if the thermal limit was not binding, wind curtailment could be significantly reduced, since the thermal plants providing inertia in Scotland could be operating simultaneously to the wind generation in that region.

In conclusion, a medium-sized fault in the low-inertia region makes it unavoidable to procure some inertia in that region, with the associated wind curtailment. 
This demonstrates the value of strong transmission corridors in systems dominated by RES generation, as the benefits are not only in terms of enhanced energy sharing across the system, but also from allowing a higher degree of freedom in the scheduling of ancillary services in each region.

\subsection{Validation of the regional-frequency-constrained SUC} \label{Sec:SimultaneousLoss}

To validate the solution of the SUC and demonstrate that it schedules sufficient ancillary services to secure the system frequencies in the event of a generation outage, two example solutions of the SUC are fed into a dynamic simulation to obtain the post-contingency frequency evolution.

\begin{figure}[!t] 
\centering 
\hspace*{-8mm}
\includegraphics[width=3.1in]{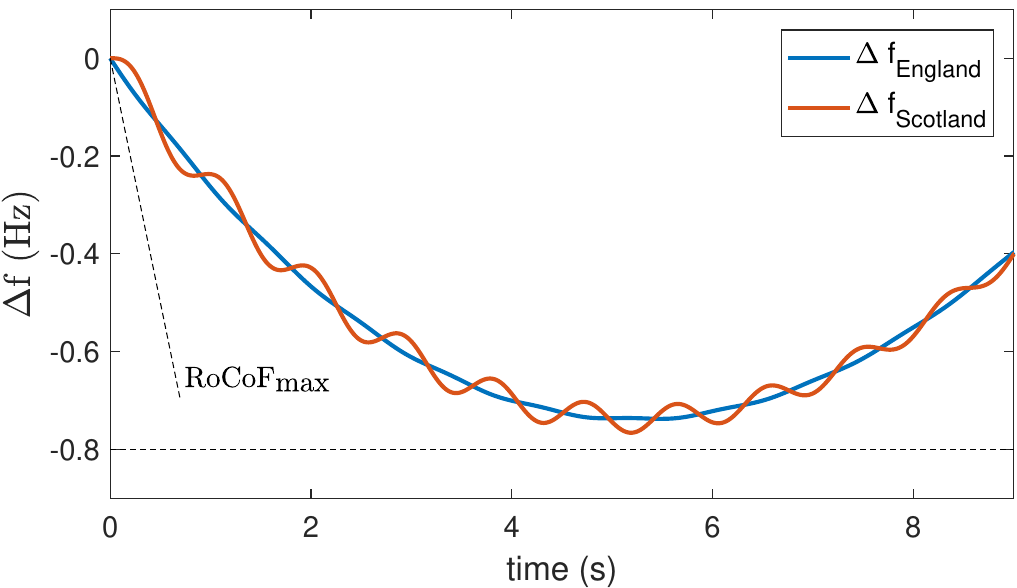}
\vspace{1mm}
\caption{Post-contingency frequency evolution, for a solution of the SUC considered in Section~\ref{Sec:LocationInertiaResponse} with a 1.8GW loss in England.}
\label{fig:DynamicSimulationSUC_FaultEngland}
\end{figure}

Fig.~\ref{fig:DynamicSimulationSUC_FaultEngland} is obtained for the case considered in Section~\ref{Sec:LocationInertiaResponse}, with $144\textrm{GW}\cdot \textrm{s}$ of inertia in England and $11\textrm{GW}\cdot \textrm{s}$ in Scotland, and 2.9GW of PFR in England and 0.25GW in Scotland. The frequency nadir is indeed above the limit in both regions, with a conservativeness of 0.03Hz (within the ranges discussed in Section~\ref{sec:NadirRegionalAssessment}), while RoCoF is not binding in this case: given the significant size of the contingency (1.8GW loss), and the fact that it takes place in the high-inertia region (which implies that inter-area oscillations are not very significant), nadir becomes the constraint driving the requirement for ancillary services.

\begin{figure}[!t] 
\centering 
\hspace*{-8mm}
\includegraphics[width=3.1in]{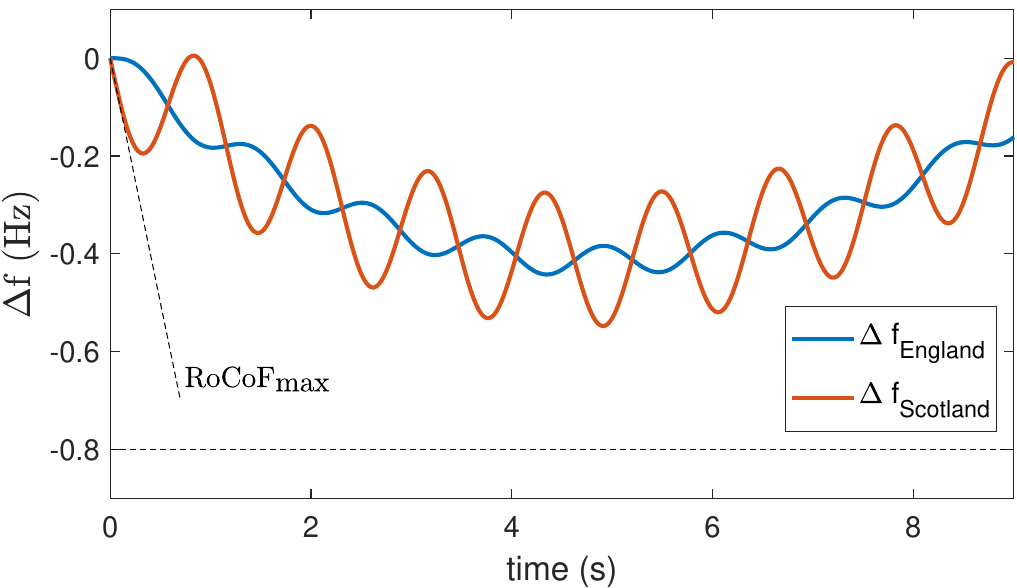}
\vspace{1mm}
\caption{Post-contingency frequency evolution, for a solution of the SUC considered in Section~\ref{Sec:FaultScotland} with a 0.8GW loss in Scotland.}
\label{fig:DynamicSimulationSUC_FaultScotland}
\end{figure}

Fig.~\ref{fig:DynamicSimulationSUC_FaultScotland} displays a case from Section~\ref{Sec:FaultScotland}, with $91\textrm{GW}\cdot \textrm{s}$ of inertia in England and $20.9\textrm{GW}\cdot \textrm{s}$ in Scotland, and 1GW of PFR in England and 0.5GW of PFR in Scotland. In this case, the amplitude of the oscillations is noticeable, which leads to RoCoF being the binding constraint in Scotland to keep this magnitude below the limit of $\textrm{RoCoF}_\textrm{max}=1\textrm{Hz/s}$. The conservativeness in RoCoF is negligible ($\sim 0.002\textrm{Hz/s}$), which is consistent with the results presented in Section~\ref{Sec:CheckConservativeRoCoF_i} since the proposed constraints introduce virtually no approximation for the RoCoF of the faulted region.

\subsection{Securing simultaneously the largest loss in both regions} \label{Sec:SimultaneousLoss}

The case studies considered in sections~\ref{Sec:LocationInertiaResponse} and \ref{Sec:FaultScotland} analyse the implications of securing against the largest loss in England or Scotland, independently for each region. However, since the system operator cannot know in advance where an outage might occur, the system scheduling must actually be secured for the fault in either region: the contingencies of 1.8GW in England and 0.8GW in Scotland must be considered simultaneously. The results for this case are presented in Table~\ref{tab:SimultaneousLoss}, showing the average values for inertia and PFR throughout the year.

The total inertia and PFR in GB are slightly above the values in Section~\ref{Sec:LocationInertiaResponse}. As inertia scheduled in Scotland is roughly the same to the one in Section~\ref{Sec:FaultScotland}, that leads to scheduling slightly less inertia in England than in Section~\ref{Sec:LocationInertiaResponse}. Given that the size of the largest loss in England is significantly higher than in Scotland, when securing for the two losses simultaneously the aggregate solution for the system is fairly similar to the one obtained for securing only the loss in England: the only difference is that a higher amount of frequency services is required in Scotland to contain the RoCoF in that region.

\begin{table}[!t]
\captionsetup{justification=centering, textfont={sc,footnotesize}, labelfont=footnotesize, labelsep=newline} 
\renewcommand{\arraystretch}{1.9}
\caption{Frequency services scheduled for securing simultaneously the largest loss in England and Scotland}
\label{tab:SimultaneousLoss}
\centering
\begin{tabular}{>{\centering\arraybackslash}m{1.9cm}|>{\centering\arraybackslash}m{1.6cm}|>{\centering\arraybackslash}m{1.6cm}|>{\centering\arraybackslash}m{1.6cm}}
    \multicolumn{1}{l|}{}  & England  & Scotland & Total GB \\ 
\hline
\raggedright Inertia ($\textrm{GW}\cdot \textrm{s}$) & 135.2 & 20.3 &  155.5 \\
\hline
\raggedright PFR (GW) & 2.9 & 0.5 & 3.4 \\
\hline

\end{tabular}
\vspace*{3mm}
\end{table}

\subsection{Computational performance of the regional-frequency-constrained SUC}

This section analyses the computational performance of the case studies presented in previous sections. The computing times presented in Table~\ref{tab:computation} were obtained using a twelve-core 3.5GHz Intel Xeon CPU with 64GB of RAM. The Stochastic Unit Commitment is formulated in a multi-threaded C++ application, while the optimisation is eventually solved by FICO Xpress 8.10 using a 0.5\% duality gap. Uncertainty in the SUC was modelled through a scenario tree with five net-demand quantiles of 0.005, 0.2, 0.5, 0.8 and 0.995 for the probability distribution of forecast error. This approach of using a relatively small number of quantiles was demonstrated by \cite{AlexEfficient} to provide a good performance for the rolling-planning SUC problem.

\begin{table}[!t]
\captionsetup{justification=centering, textfont={sc,footnotesize}, labelfont=footnotesize, labelsep=newline} 
\renewcommand{\arraystretch}{1.9}
\caption{Computation time for one year of system operation}
\label{tab:computation}
\centering
\begin{tabular}{>{\centering\arraybackslash}m{1.6cm}|>{\centering\arraybackslash}m{1.6cm}|>{\centering\arraybackslash}m{1.6cm}|>{\centering\arraybackslash}m{1.6cm}}
    Case in Section \ref{Sec:LocationInertiaResponse}  & Case in Section \ref{Sec:FaultScotland}  & Case in Section \ref{Sec:SimultaneousLoss} & Extra case, $\quad$4 Regions\\ 
\hline
22.5h & 15min & 21.5h & 50min \\
\hline

\end{tabular}
\end{table}

The case in Section \ref{Sec:LocationInertiaResponse} considers the two-region system England-Scotland, secured against a loss of 1.8GW in England, while the case in Section \ref{Sec:FaultScotland} considers the same two-region system but secured for a contingency of 0.8GW in Scotland. The latter is solved in a significantly shorter time than the former, due to the binding constraints in each case: RoCoF constraints are the most challenging to meet for the 0.8GW loss in Scotland, while nadir constraints drive the requirement for inertia and PFR in the case of 1.8GW loss in England. 
Since the nadir constraints include binary variables for the discretisation, this case has a higher computation time, as MILPs with increasing number of binary variables are typically harder to solve. 

An extra case with four regions has also been considered, where all four regions have the same characteristics as England in terms of generation mix, demand and size of largest possible contingency (i.e.~1.8GW). The SUC is secured for the largest loss in each of the regions. The computation time for this case is under 1h, due to the fact that inertia is needed in every region to secure the local 1.8GW loss, which implies a high total system inertia that makes the nadir constraints non-binding. This demonstrates that systems with more regions for frequency purposes (i.e.~systems where machines are split in a higher number of coherent groups) do not necessarily increase computation time: the burden of the optimisation shows to be mainly driven by the constraints that are binding in a given system, that is, RoCoF or nadir. 

The results in Table~\ref{tab:computation} demonstrate that the already computationally-expensive stochastic scheduling is still tractable when including the regional frequency constraints. The formulation of this SUC is an MILP, which can be handled by the vast majority of commercial as well as open-source solvers. For a particular system configuration in which the nadir constraints are known to be binding, decomposition techniques for the optimisation problem could be explored to speed up computation (e.g.~Benders decomposition, ADMM), which are however outside of the scope of this paper.

\section{Conclusion and Future work} \label{Sec:Conclusion}

This paper has demonstrated the applicability of the frequency-stability conditions proposed in Part I to be implemented in optimisation problems formulated as MILPs. The Great Britain system has been used as the platform to test the proposed frequency-secured framework, demonstrating that it is key to procure inertia and response appropriately among the different regions of the system, and that medium-sized faults can have significant impacts if they occur in low-inertia regions. While the quantitative results presented in this paper apply to the GB 2030 system that has been considered, the qualitative findings summarised here apply to any system with high penetration of non-synchronous generation located in isolated or weakly interconnected regions. It is credible that, to a certain degree, that will be the case for several other systems in the world as they become decarbonised to meet emissions targets.

Future lines of work on this topic should explore how to include probabilistic metrics for the provision of frequency response (so as to consider in a realistic way the response from distributed providers) and how to simultaneously guarantee transient stability in low-inertia systems after a generation outage (since the post-fault activation of response must not cause a phase-angle separation that could potentially lead to loss of synchronism).


\ifCLASSOPTIONcaptionsoff
  \newpage
\fi

\IEEEtriggeratref{30}

\bibliographystyle{IEEEtran} 
\bibliography{Luis_PhD}

\vskip -1\baselineskip plus -2fil

\begin{IEEEbiographynophoto}{Luis Badesa}
(S'14-M'20) received the Ph.D. degree in Electrical Engineering from Imperial College London, U.K., in 2020. He is currently a Research Associate within the Control \& Power research group at Imperial College London. His research interests lie in modelling and optimisation for low-carbon power grids.
\end{IEEEbiographynophoto}

\vskip -1\baselineskip plus -2fil

\begin{IEEEbiographynophoto}{Fei Teng}
(M'15) received the Ph.D. degree in Electrical Engineering from Imperial College London, U.K., in 2015. Currently, he is a Lecturer in the Department of Electrical and Electronic Engineering, Imperial College London, U.K. His research focuses on scheduling and market design for low-inertia power system, cyber-resilient energy system operation and control, and objective-based data analytics for future energy systems.
\end{IEEEbiographynophoto}

\vskip -1\baselineskip plus -2fil

\begin{IEEEbiographynophoto}{Goran Strbac}
(M'95) is Professor of Electrical Energy Systems at Imperial College London, U.K. His current research is focused on the optimization of operation and investment of low-carbon energy systems, energy infrastructure reliability and future energy markets. 
\end{IEEEbiographynophoto}

\end{document}